\documentclass[12pt]{article}
\usepackage{pdfsync}
\usepackage[utf8]{inputenc}
\usepackage[T1]{fontenc}
\usepackage{lmodern}
\usepackage{amsmath,amsthm,amsfonts}
\usepackage{enumitem}

\newtheorem{theorem}{Theorem}
\theoremstyle{remark}
\newtheorem{remark}{Remark}
\newtheorem{prop}{Proposition}
\def\xN{\mathbb{N}}
\def\xZ{\mathbb{Z}}
\newcommand{\abs}[1]{\left|#1\right|}
\newcommand{\nrm}[2]{\left\|#1\right\|_{#2}}

\def\xR{\mathbb{R}}
\newcommand{\nbOne}{\mathbf{1}}        
\newcommand{\ind}[1]{\nbOne_{#1}}  

\newcommand{\overbar}[1]{\mkern 1.5mu\overline{\mkern-1.5mu#1\mkern-1.5mu}\mkern 1.5mu}

\newcommand{\e}{\epsilon}

\newcommand{\ud}{\,d}

\newcommand{\prbPt}[2]{\mathbb{P}_{#1}\left[#2\right]}
\newcommand{\prb}[1]{\mathbb{P}\left[#1\right]}       
\newcommand{\esp}[1]{\mathbb{E}\left[#1\right]}       

\newcommand{\NRM}[1]{{{\left\| #1\right\|}}} 
\newcommand{\TVNRM}[1]{\NRM{#1}_{\mathrm{TV}}} 
\newcommand{\PAR}[1]{{{\left(#1\right)}}} 

\begin{document}

\title{Piecewise deterministic Markov process --- recent results}
\author{Romain Aza\"{i}s\thanks{Inria Bordeaux Sud-Ouest, team CQFD et Universit\'e de Bordeaux, IMB, CNRS
UMR 5251, 200, Avenue de la Vieille Tour, 33405 Talence cedex, France.}\and Jean-Baptiste Bardet\thanks{Laboratoire de Math\'ematiques Rapha\"el Salem, Universit\'e de Rouen, Avenue de l'Universit\'e, BP 12, 76801 Saint-Etienne-du-Rouvray, France.}\and Alexandre G\'enadot\thanks{Laboratoire de Probabilit\'es et mod\`eles al\'eatoires, Universit\'e Pierre et Marie Curie, 4, place Jussieu, 75005 Paris, France.}\and Nathalie Krell\thanks{Universit\'e de Rennes 1,
 Institut de Recherche math\'ematique de Rennes,
CNRS-UMR 6625, Campus de Beaulieu.
 B\^atiment 22, 35042 Rennes Cedex, France.}\and Pierre-Andr\'e Zitt\thanks{Laboratoire d'Analyse et de Math\'ematiques Appliqu\'ees,
  Universit\'e de Marne-la-Vall\'ee,
  5, boulevard Descartes,
  Cit\'e Descartes - Champs-sur-Marne,
77454 Marne-la-Vall\'ee Cedex 2}}

\maketitle

\begin{abstract}
We give a short overview of recent results on a specific class of
Markov process: the Piecewise Deterministic Markov Processes (PDMPs).
We first recall the definition of these processes and give some general results.
On more specific cases such as the TCP model or a model of switched vector
fields, better results can be proved, especially as regards long time behaviour.
We continue our review with an infinite dimensional example
of neuronal activity. From the statistical point of view,
these models provide specific challenges: we illustrate this point
with the example of the estimation of the distribution of the inter-jumping times.
We conclude with a short overview on numerical methods used for simulating PDMPs.
\end{abstract}


\section{General introduction}
The piecewise deterministic Markov processes (denoted PDMPs) were
first introduced in the literature by Davis (\cite{Dav84,Dav93}).
Already at this time, the theory of diffusions had such
powerful tools as the theory of It\=o calculus and stochastic
differential equations at its disposal. Davis's goal was to endow the
PDMP with rather general tools. The main reason for that was to
provide a general framework, since up to then only very particular
cases had been dealt with, which turned out not to be easily
generalizable.

PDMPs form a family of c\`adl\`ag Markov processes involving a
deterministic motion punctuated by random jumps. The motion of the
PDMP~$\{X(t)\}_{t\geq 0}$ depends on three local characteristics, namely
the jump rate $\lambda$, the flow $\phi$ and the transition measure Q
according to which the location of the process at the jump time is
chosen.  The process starts from $x$ and follows the flow $\phi(x, t)$
until the first jump time $T_{1}$ which occurs either spontaneously in
a Poisson-like fashion with rate $\lambda(\phi (x,t))$ or when the
flow $\phi(x, t)$ hits the boundary of the state-space. In both cases,
the location of the process at the jump time $T_{1}$, denoted by
$Z_1=X(T_{1})$, is selected by the transition measure $Q(\phi
(x,T_{1}), \cdot) $ and the motion restarts from this new point as
before. This fully describes a piecewise continuous trajectory for
$\{X(t)\}$ with jump times $\{T_{k}\}$ and post jump locations
$\{Z_{k}\}$, and which evolves according to the flow $\phi$ between
two jumps.

These processes have been heavily studied both from a theoretical and from an
applied perspective in various domains such as communication networks with the
control of congestion TCP/IP (V.~Dumas and al. \cite{DGR02}, V.~Guillemin and
al. \cite{GRZ04}), neurobiology for the Hodgkin-Huxley model of  neuronal
activity (K.~Pakdaman and al. \cite{PTW10}), reliability (F.~Dufour and
Y.~Dutuit \cite{DD02}), biologic population models  (H.G.~Othmer and
al.\cite{ODA88} as well as R.~Erban and H.G.~Othmer \cite{EO04}),
to present just a few examples.

The paper is organized as follows. In Section~\ref{sec_definition} we give a
precise definition and some general properties of the PDMPs. Then we illustrate
the state of the art regarding PDMPs through three specific examples: a model of
switched vector fields (Section~\ref{sec_switch}), the TCP process
(Section~\ref{sec_tcp}), and a modelization of neuronal activity
(Section~\ref{sec_neuron}). Finally, we briefly review some results about a
non-parametric statistical method to get an estimate of the conditional density
associated with the jumps of a PDMP defined on a separable metric space
(Section~\ref{sec_estimation}) and we end with a survey of numerical methods in
Section~\ref{sec_numerics}.

\section{Definition and some properties of PDMPs}
\label{sec_definition}

\subsection{Definition of a PDMP.}
Let $M$ be an open subset of $\xR^{n}$, $\partial M$ its boundary,
$\overbar M$ its closure and $\mathcal{B}(M)$ the set of real-valued,
bounded, measurable functions defined on~$M$. A PDMP is determined by
its local characteristics $(\phi, \lambda ,Q)$ where:
\begin{itemize}
\item The flow $\phi:\xR^{n}\times\xR\rightarrow\xR^{n}$ is a
one-parameter group of homeomorphisms: $\phi$ is continuous, $\phi (\cdot,
t)$ is an homeomorphism for each $t\in\xR$, satisfying the semigroup
property: $\phi(\cdot, t+ s) = \phi (\phi(\cdot, s), t)$.

For each $x$ in $M$, we introduce the deterministic hitting time of
the boundary
\begin{equation}
  \label{eq=defTStar}
t^{\star}(x):=\inf\{t > 0 :\phi(x, t) \in\partial M\},
\end{equation}
with the convention $\inf \emptyset =\infty$.

\item The jump rate $\lambda :\overbar{M} \rightarrow \xR_{+}$ is assumed
to be a measurable function satisfying
\[
\forall x\in M,\; \exists\,\epsilon>0\:\; \text{such that}\:\;
\int_{0}^{\epsilon}\lambda(\phi (x,s))ds<\infty.
\]

\item $Q$ is a Markov kernel on $(\overbar M,\mathcal{B}(\overbar M))$
satisfying the following property:
\[
\forall x\in \overbar M,\quad Q(x,M-\{x\})=1.
\]
\end{itemize}
From these characteristics, it can be shown \cite[pp 62--66]{Dav93}
that there exists a filtered probability space $(\Omega ,
\mathcal{F},\{\mathcal{F}_{t}\}, \{\mathbb{P}_{x}\})$ such that the
motion of the process $\{X(t)\}$ starting from a point $x\in M$ may be
constructed as follows.  Consider a random variable $T_{1}$ such that
\[
\mathbb{P}_{x}\{T_{1}>t\}=
  \begin{cases}
   e^{-\Lambda (x,t)}  &\text{ for $t<t^{\star}(x)$,}\\
   0                   &\text{ for $t\geq t^{\star}(x)$,}
 \end{cases}
\]
where for $x\in M$ and $t \in [0, t^{\star}(x)]$
\[
\Lambda (x,t)=\int_{0}^{t}\lambda (\phi (x,s))ds.
\]
If $T_{1}$ is equal to infinity, then the process $X$ follows the
flow, i.e.  $t\in \xR_{+}$, $X(t) = \phi (x, t)$. Otherwise select
independently a $M$-valued random variable (labelled $Z_{1}$) having
distribution $Q(\phi (x,T_{1}),\cdot)$, namely $\mathbb{P}_{x}(Z_{1}\in
A) = Q(\phi (x,T_{1}),A)$ for any $A \in\mathcal{B}(M)$. The
trajectory of $\{X(t)\}$ starting at $x$, for $t \in [0,T_{1}]$, is
given by
\[
X(t)=
\begin{cases}
\phi (x,t) &\text{for } t<T_{1},\\
Z_{1}      &\text{for } t= T_{1}.
\end{cases}
\]
Starting from $X(T_{1}) = Z_{1}$, we now select the next inter-jump
time $T_{2}- T_{1}$ and post-jump location $X(T_{2}) = Z_{2}$ in a
similar way.

This construction properly defines a Markov process $\{X(t)\}$ which
satisfies the strong Markov property with jump times
$\{T_{k}\}_{k\in\mathbb{N}}$ (where $T_{0} = 0$).  A very natural
Markov chain is linked to $\{X(t)\}$, namely the chain
$(\Theta_{n})_{n\in\mathbb{N}}$ defined by $\Theta_{n} =
(Z_{n},S_{n})$ with $Z_{n} =X(T_{n})$ and $S_{n} = T_{n}-T_{n-1}$ for
$n\geq 1$ and $S_{0} = 0$. Clearly, the process
$(\Theta_{n})_{n\in\mathbb{N}}$ is a Markov chain.  This chain will
turn out to be particularly useful in the next sections (see
Section~\ref{sec_numerics}).

\begin{remark}
  \label{rmrk:espace_etat}
Davis originally defined  PDMPs on a
disjoint union  $\bigcup_{v\in K} \{v\}\times M_v$,
where $K$ is a countable index set and, for each $v$,
$M_v$ is a subset of $\xR^{n}$.  The
definition above, with a single copy of $\xR^n$, can be used without
loss of generality -- see Remark 24.9 in \cite{Dav93} for details.
Depending on the process, one or the other definition may be more
natural: we will use the original definition in Section
\ref{sec_switch}.
\end{remark}

\subsection{A useful discrete process}
A basic assumption in all the following will be that,
for every starting point $x\in M$,  $\mathbb{E}
\left\{\sum_{k}\nbOne_{\{t\geq T_{k}\}}\right\}<\infty$.
This ensures the non explosion of the process.

The chain $\Theta_n$ obtained by observing the process at jump
times may not have enough jumps to guarantee good comparison
properties with the continuous time process~$Z_t$. To address this
problem, O.~Costa and F~Dufour introduced in~\cite{CD99} another discrete chain
as follows. The idea is to record the positions of the continuous process
both at jump times and at additional random times given by an
independent Poisson process of rate $1$. Formally,
first define two  substochastic kernels $H$ and $J$ by
\begin{align*}
H(x,A)& := \int_{0}^{t_{\star}(x)} e^{-(s+\Lambda(x,s))} \ind{A}(\phi(x,s)) \ud s,\\
J(x,A)&:=\int_{0}^{t_{\star}(x)}\lambda (\phi(x,s))e^{-(s+\Lambda(x,s))}Q(\phi(x,s),A)ds
+e^{-(t_{\star}(x)+\Lambda(x,t_\star(x)))}Q(\phi(x,t_{\star}(x)),A).
\end{align*}
The first kernel corresponds to the additional observation times,
the second one to the "real jumps", either due to the jump rate
(first term) or the hitting of the boundary (second term).
The sum $G = J+H$ is proved in~\cite{CD99} to be a Markov kernel.
We denote by $\{\theta_{n}\}$ the associated Markov chain,
which can be generated from the sample paths of $\{X_t\}$
by adding observation times as said before -- see~\cite[Theorem 3.1]{CD08}
for a formal statement.

\begin{theorem}[Discrete and continuous processes, \cite{CD08}]
  \label{th:CostaEtDufour}%
  The Markov chain $\{\theta_{n}\}$ and the original process $\{X_t\}$ are closely
related:
\begin{enumerate}
\item The PDMP $\{X(t)\}$ is irreducible if and only if the Markov
chain $\{\theta_{n}\}$ is irreducible.
\item If $\nu$ is an invariant measure for
$\{X(t)\}$, then $\nu\sum_{0}^{\infty} J^{j}$ is invariant for
$\{\theta_{n}\}$ and $\nu\sum_{0}^{\infty} J^{j}=\nu$. Conversely,
if~$\pi$ is an invariant measure for $\{\theta_{n}\}$, then $\pi
H$ is invariant for $\{X(t)\}$ and $\pi\sum_{0}^{\infty}
J^{j}=\pi$.
\item The PDMP $\{X(t)\}$ is recurrent if and only if the Markov chain
$\{\theta_{n}\}$  is recurrent.
\item The PDMP $\{X(t)\}$ is positive recurrent if and only if the
Markov chain $\{\theta_{n}\}$  is recurrent with
invariant measure $\pi$ satisfying $\pi H(M) < \infty$.
\item The PDMP $\{X(t)\}$ is Harris recurrent if and only if the
Markov chain $\{\theta_{n}\}$ is Harris recurrent.
\end{enumerate}
\end{theorem}

The authors in~\cite{CD08} also give sufficient conditions on~$G$ (in
a modified Foster-Lyapunov criterion form) to ensure the existence of
an invariant probability measure, positive Harris recurrence and
ergodicity for the PDMP.

In \cite{Las04}, stability and ergodicity via Meyn-Tweedie arguments
are established for one dimensional PDMPs as AIMD (Additive Increase
Multiplicative Decrease) processes. Most of these results are
essentially of qualitative type, i.e. no practical information can be
obtained on the rate of convergence to equilibrium. The papers
\cite{CMP10,FGM12} are first attempts to get quantitative
results for the long time behavior of special PDMPs: TCP (which is
defined in Section~\ref{sec_tcp}) for the first one and a PDMP
describing the motion of a bacteria for the second one.

\section{Regularity and general convergence results for Markov switching systems}
\label{sec_switch}

\subsection{The Markov switching model}
We consider a subclass of PDMP sometimes known as ``Markov switching model''.
The ``natural'' state space of this process is the product space $\xR^d
\times\{1, \ldots, n\}$ --- as was said in Remark~\ref{rmrk:espace_etat}, this
fits within the general theory described above.  To define the process,
consider $n$ vector fields $F^1$, $F^2$, \ldots, $F^n$ on $\xR^d$. These fields
define flows $\phi^i$ via the ODE:
\[
\frac{dy_t}{dt} = F^i(y_t) dt.
\]
To simplify matters we suppose that the $F^i$ are $\mathcal{C}^1$
and that there is a compact set $K\subset \xR^d$ that is left invariant by all flows
(nasty things may occur if this is not the case, see e.g.{} \cite{BLMZ1}).

We also suppose that we are given $n^2$ nonnegative functions
$\lambda_{ij}:\xR^d \to \xR$ (the jump rates), such that
$\lambda_{ii}(y) = 0$, and for any given $y$,
$(\lambda_{ij}(y))_{i,j}$ is irreducible.

The Markov switching model is a process $X_t = (Y_t,I_t) \in
\xR^d\times \{1,2,\cdots n\}$ defined informally as follows:
\begin{itemize}
  \item $Y_t$ is driven by the vector field $F^{I_t}$,
  \item if $I_t = i$ and $Y_t = y$, $I_t$ jumps to $j$ with rate
$\lambda_{i,j}(y)$.
\end{itemize}
The corresponding generator is:
\[
  Lf(y,i) = F^i \cdot \nabla_x f + \sum_j \lambda_{ij}(y) (f(y,j) -
f(y,i)).
\]
This process can be constructed as described in
Section~\ref{sec_definition}.  Another possible way is to generate a
Poisson process with a high enough intensity to ``propose'' jump
times, and to accept or reject these jumps with probabilities
depending only on the values of $\lambda_{ij}(y)$ at the current point
(see Section~2 of~\cite{BLMZ2} for details). This naturally defines a
discrete process $\tilde{X}_n$ (a Markov chain), which turns out to be
an alternative to the chain $\theta_n$ defined above, in the sense
that results similar to  Theorem~\ref{th:CostaEtDufour} hold.

In this model, the trajectory of the ``position'' $Y_t\in\xR^d$ does not jump,
it evolves continuously. The jumps only occur on the discrete part $I_t$, the
one that dictates which (deterministic) dynamics the position must follow. Let
us also note that there are no ``boundary jumps'' in this model: for every
$x=(y,i)$, the quantity $t^\star(x)$ defined by~\eqref{eq=defTStar} is infinite.

Many questions that seem intractable for fully general PDMPs can be answered in
this restricted framework. In particular, a kind of regularity for the law of
$(Y_t,I_t)$ can be established if the vector fields satisfy a condition that
closely resembles Hörmander's classical condition for diffusion. These results
are described in the following section. In the next one we discuss results on
the ``infinite time'' behaviour of the process: once more, exponential
convergence to equilibrium can be proved under fairly weak assumptions.

\subsection{Regularity results}
We are interested here in criteria that guarantee
a regularity for the law of the continuous part, $Y_t$.

For the sake of comparison, let $B_t= (B^1_t, \ldots, B^n_t)$ be an
$n$-dimensional Brownian motion, and consider the solution of the
following (Stratonovich) SDE:
\[
  dY_t = F^0(Y_t)dt + \sum_{i=1}^n F^i(Y_t)\circ dB^i_t.
\]
Depending
on the $F^i$, the law of $Y_t$ may or may not have a density with
respect to Lebesgue measure. For example, if the $F^i$ are constant,
the law is regular if and only if the $F^i$ span $\xR^d$.  For
non-constant fields there is a well-known criterion for regularity,
known as Hörmander's criterion, and expressed in terms of the $F^i$
and their iterated Lie brackets. Thanks to this criterion, one may
prove that the law of a kinetic process $(Y,V)$ driven only by a
Brownian motion on its velocity is immediately absolutely continuous,
with a regular density (see \textit{e.g.} \cite{villahypo}).

In our case, there is no hope for such a strong result.  Indeed, the
probability that $(Y_t,I_t)$ has not jumped at all between $0$ and $t$
is strictly positive, therefore the law of $Y_t$ has a singular part
-- we will see a similar phenomenon in the study of the TCP process
with the lower bound \eqref{eq:lowerTV}.

However, let us define the following families of vector fields:
\begin{align*}
  \mathcal{V}_0 &= \{F^i - F^j, i\neq j\}, \\
  \forall i\geq 0, \qquad \mathcal{V}_{i+1} &
= \mathcal{V}_i \cup \{ [F^i, G], G \in \mathcal{V}_i\}.
\end{align*}
With these notations the following result holds:
\begin{theorem}
  [Absolute continuity, \cite{BH,BLMZ2}]
  \label{thrm:regularity}
  Suppose there is a point $p\in\xR^d$ and an integer $k$ such that
the ``bracket condition'' holds, i.e.\ the iterated Lie brackets
$\{V(p), V \in \mathcal{V}_k\}$ span $\xR^d$.  Then there is a time
$t>0$, a strictly positive $\eta$, and two probability measures
$\nu^c_t, \nu^r_t$ on $\xR^d\times\{1,\ldots,n\}$  such that:
  \begin{itemize}
  	\item the law $\mu_t$ of $Y_t$ starting from $p$ can be
written as $\mu_t = \eta \nu^c_t + (1-\eta)\nu^r_t$,
	\item the ``continuous'' part $\nu^c_t$ is absolutely
continuous with respect to the product of the Lebesgue measure
on $\xR^d$ and the counting measure on $\{1,\dots,n\}$.
  \end{itemize}
\end{theorem} Similar results are obtained in \cite{BH,BLMZ2} for the
discrete chain $\tilde{Z}_n$; they involve another (weaker) bracket
condition.

\subsection{Invariant measures and convergence}
To study the long time behaviour of $(Y_t,I_t)$, it is natural
to look for an analogue of the $\omega$-limit set (defined when there
is only one flow). To that end,
first define the positive trajectory  of $y$ as the set
\[
\gamma^+(y) =
\bigcup_{m\in\xN} \left\{
    \phi^{i_m}_{t_m}(\cdots (\phi^{i_1}{t_1}(y))\cdot )
    \middle|
    (i_1,\ldots, i_m) \in \{1, \ldots, n\}^m;
    (t_1, \ldots, t_m) \in \xR_+^m
  \right\}.
\]
That is, $\gamma^+(y)$ is the set of points reachable from $y$ by a
trajectory of the process. Recall that all the $F^i$ leave
a compact set $K$ invariant. The right limit set to consider is
the \emph{accessible set},  the  (possibly empty) compact
set $\Gamma \subset K$ defined as
\[
\Gamma = \bigcap_{x \in K} \overline{\gamma^+(x)}.
\]
Informally this is the set of points that are (approximately)
reachable from anywhere in $K$ by following the flows.

This set is deeply connected with the invariant measures of the
process.  Let $\mathcal{P}_{inv}$ be the set of probability measures
that are invariant for the process, and let $E=\{1, \ldots, n\}$ be
the finite index set.
\begin{prop}[Limit set and support of invariant measures]
  \ 
\begin{enumerate}[label={\bfseries(\roman*)}]
  \item If $\Gamma \neq \emptyset$ then $\Gamma \times E\subset
supp(\mu)$ for all $\mu\in\mathcal{P}_{inv}$ and there exists $\mu \in
\mathcal{P}_{inv}$ such that $supp(\mu) = \Gamma \times E.$
\item If $\Gamma$ has nonempty interior, then $\Gamma \times E =
supp(\mu)$ for all $\mu \in \mathcal{P}_{inv}$.
\item Suppose that $\mathcal{P}_{inv}$ is a singleton $\{\pi\}$.  Then
$supp(\pi) = \Gamma \times E.$
\end{enumerate}
\end{prop}
Once more, similar statements hold for the invariant
measures of the discrete chain $\tilde{Z}_n$.

A usual way of proving convergence to equilibrium for Markov chains or
processes, following the classical idea of the Foster-Lyapunov
criterion, is to identify a ``good'' subset $A$ of the state space,
satisfying two properties:
\begin{itemize}
	\item the process has a tendency to return to $A$;
	\item for any $x$, $y$ in $A$, the laws of the process
	  starting from $x$ and $y$ are similar in some sense.
\end{itemize}
If these two properties hold, one can construct a coupling between
trajectories in three steps:
\begin{itemize}
	\item wait for the two processes to come back to $A$;
	\item once they are in $A$, using the second property to try
	  to couple them;
	\item if the coupling fails, go back to step $1$.
\end{itemize}
In the classical case of a discrete Markov chain, $A$ will usually be
a single point so that the second property will be trivial; the
``tendency to return'' will typically follow from the existence of a
Lyapunov function. This type of technique is refined and generalized
for Markov chains on general state spaces in the renowned book by Meyn
and Tweedie~\cite{MT93}.

In our case, without additional assumptions, we cannot hope to find a
Lyapunov function. However, if $\Gamma$ is not empty, and if $p$ is a
point in $\Gamma$, we can prove --- using compactness arguments ---
that the process will visit infinitely often any neighborhood of $p$.

To use the heuristics above, it remains to show that two copies of the
process, starting from a point near $p$, have a good chance to
couple. This is where the regularity results may be used: two
variables with strictly positive density can be coupled with positive
probability.  This is the main idea of the proof of:
\begin{theorem}
  [Exponential convergence, \cite{BLMZ2}]
  Suppose that the ``bracket condition'' of
Theorem~\ref{thrm:regularity} is satisfied at a point $p\in\Gamma$.
Then there exists a unique invariant probability measure $\mu$, its
support is $\Gamma\times \{1,\ldots, n\}$, and $X_t = (Y_t,I_t)$
converges exponentially fast in total variation norm:
  \begin{equation}
    \exists\, \alpha>0, C>0, \quad \nrm{\mathcal{L}(X_t) - \pi}{\mathrm{TV}}
\leq C \exp(-\alpha t).
    \label{eq=convergenceContinue}
  \end{equation}
\end{theorem}

These convergence results are quite general. However, there is no hope
to find reasonable convergence rates by this approach, since we used
compactness arguments.  In particular cases, one can often use the
specific structure of the process to prove much better bounds. This is
done in the Markov switching case under additional assumptions in
\cite{BLMZ12}. The next section describes similar results for the
slightly different case of the TCP window-size process.

\section{Explicit convergence rates via coupling}
\label{sec_tcp}

Beyond the problems of existence, uniqueness, and regularity of an
invariant probability measure, an important and challenging question
is to get explicit and efficient rates of convergence to equilibrium
for PDMPs.  It seems hard to address the problem via spectral methods
since a PDMP is generally, and inherently, non reversible (see however
\cite{GM13,MM12} for positive results in this direction). As mentioned
previously, coupling methods ``\`a la'' Meyn-Tweedie can be applied but
hardly give efficient rates of convergence. We will present here a
specific method that allows to treat one specific (and relatively
simple) example of PDMP, the so-called TCP window-size process.

\subsection{The TCP window-size process}

The TCP window-size process appears as a scaling limit of the
transmission rate of a server uploading packets on the Internet
according to the algorithm used in the TCP (Transmission Control
Protocol) in order to avoid congestions (see~\cite{DGR02} for details
on this scaling limit). This process $X=\{X_t\}_{t\geq0}$ has
$[0,\infty)$ as state space and its infinitesimal generator is given,
for any smooth function $f:[0,\infty)\to\xR$, by
\begin{equation}
\label{eq:G1}
  Lf(x) = f'(x) + x(f(x/2)-f(x)).
\end{equation}
This means that the transmission rate is increased linearly when there
is no congestion, and divided by two in case of congestion. A
congestion occurs at random times, with rate given by the value of the
transmission rate. This process is rather simple and its stationary
regime is well understood: it is ergodic and admits a unique invariant
probability measure $\mu$, with an explicit density (see~\cite{DGR02}).
Getting an explicit speed of convergence is however not so
simple, since the TCP process is irreversible but also because, in
some sense, it ``lacks randomness'' (the only randomness comes from
the jump times), what makes any coupling procedure delicate,
particularly for convergence in total variation norm.  In the
following, the semigroup associated to the TCP process will be denoted
$\{P_t\}_{t\geq0}$.

\subsection{Speed of convergence in Wasserstein distance}

The speed of convergence in Wasserstein distance can be dealt with
using a coupling first introduced in \cite{CMP10}. This is a Markov
process on $\xR_+^2$ whose marginals are two TCP processes, and which
is defined by its generator as follows: if $x\geq y$,
\begin{equation}
\label{eq:wascou}
\mathfrak{L} f(x,y)
=(\partial_{x}+\partial_{y})f(x,y)+y\big(f(x/2,y/2)-f(x,y)\big)
+(x-y)\big(f(x/2,y)-f(x,y)\big),
\end{equation}
and it is given by the symmetric expression if $x>y$. In words, with a
rate equal to the minimum of the two coordinates, both trajectories
jump simultaneously, and with a rate equal to the difference between
the two coordinates the higher trajectory jumps alone (and the lower
one never jumps alone). It induces problems since the higher
trajectory may jump alone even when the distance between both
trajectories is small (hence at a small but positive rate), increasing
dramatically this distance. This can also be seen on computation,
noting that this dynamics does not contract the distance between
trajectories in $L^p$ distance, for any $p\geq1$. However, one gets
(non uniform) contraction for $p=1/2$. Working this out, and using
also a bound on any moment of the process at any time uniformly in the
starting point, one finally gets the following result:

\begin{theorem}
[Explicit rate for TCP in Wasserstein distance, \cite{BCGMZ}]
\label{th:wasserstein}
Let  $c=\sqrt{2}(3+\sqrt{3})/8 \approx 0.84$
and
$\lambda=\sqrt 2(1-\sqrt c)\approx 0.12$.
For any $\tilde \lambda < \lambda$, any  $p\geq1$ and any $t_0>0$,
there is a constant $C=C(p,\tilde \lambda,t_0)$ such that,
for any initial probability measures $\nu$ and $\tilde \nu$ and
any $t\geq t_0$,
\[
W_p(\nu P_t,\tilde \nu P_t)\leq C \exp\left(-\frac{\tilde \lambda}pt\right)
.
\]
\end{theorem}

This result can be compared to numerical simulations: the exponential
rate of convergence for the $W_1$ Wasserstein distance to the
invariant measure is about $1.6$, and the one for the $L^1$ distance
between both coordinates of the coupling we used is about $0.5$ (which
means that one cannot do better using this coupling).

\subsection{Speed of convergence in total variation norm}

Let us first note that, even starting from two arbitrarily close
initial conditions, there is a positive probability that two TCP
processes never have jumped before a given time $t$, which gives an
immediate and striking lower bound for the total variation distance
between laws at time $t$. To be more specific, the first jump time
$T_1$ of the TCP process starting from $x$ is characterized by
$p_t(x)=\mathbb{P}\{T_1\geq t\}=e^{-t^2/2-xt}$. This implies that
\begin{equation}
\label{eq:lowerTV}
\|\delta_xP_t-\delta_yP_t\|_{\textrm{TV}}\geq
p_t(x)\vee p_t(y)=e^{-t^2/2-(x\wedge y) t}\,.
\end{equation}
Consequently, this distance is not tending to 0 as $y$ tends to
$x$. Nevertheless, it can be shown that, for $y$ close to $x$, the
distance is really on the order of this lower bound. To note that, we
construct a coupling which tries to stick two TCP processes in one
jump. Indeed, assume that $x>y$ and denote by $T_1^x$ the first jump
time of the process starting from $x$ and $T_1^y, T_2^y$ the two first
jump times of a TCP process starting from $y$. Then if these random
variables satisfy the two relations
\begin{equation}
\label{eq:coupltps}
T_1^x=T_1^y+x-y
\end{equation}
and $T_2^y-T_1^y\geq x-y$, both processes are at the same position at
time $T_1^x$, hence can be coupled from this time. The higher
probability of getting \eqref{eq:coupltps} is an optimal coupling
probability for two measures with density, it can hence be computed
explicitly, and gives an upper bound on the order of
\eqref{eq:lowerTV}.

We finally can construct a coupling in the following way: for any
$\varepsilon>0$, we use the dynamical coupling defined in
\eqref{eq:wascou} during a time $t_1$ on the order of
$\log(1/\varepsilon)$ to try to make the trajectories $\varepsilon$
close; we then try to stick them in a time $t_2$ on the order of
$\sqrt{\log(1/\varepsilon)}$. Taking $\varepsilon $ to $0$ gives the
following result:

\begin{theorem}
[Explicit rate for TCP in total variation, \cite{BCGMZ}]
\label{th:variation}
Let $\lambda$ as defined in Theorem \ref{th:wasserstein}. For any
$\tilde \lambda<\lambda$ and any $t_0>0$, there exists $C$ such that,
for any initial probability measures $\nu$ and $\tilde \nu$ and any
$t\geq t_0$,
\[
\TVNRM{ \nu P_t-\tilde \nu P_t}\leq
C\exp\PAR{-\frac{2\tilde \lambda}{3}t}\,.
\]
 \end{theorem}

\subsection{Application to other processes}
This strategy (first apply a coupling which is efficient for the
Wasserstein distance, then stick the trajectories in one single
attempt) can also be used for other processes. We can treat for
example the case of the TCP process with constant jump rate, with
generator
\begin{equation}
\label{eq:G2}
  Lf(x) = f'(x) + r(f(x/2)-f(x)),
\end{equation}
$r>0$ being fixed. This gives a purely probabilistic proof (and
slightly better constants) of results first obtained in an analytic
way by Perthame and Ryzhik in \cite{PR}.

The same idea can also be applied to diffusion processes, associating
recent results from \cite{eberle,BGG} on the convergence in
Wasserstein distance with results on the regularity of the semigroup
at small time from \cite{Wang10}.  The reader is referred to
\cite{BCGMZ} for details.

\section{Infinite dimensional PDMP}
\label{sec_neuron}

\subsection{Spatially extended conductance based neuron model}
The intrinsic variability of a large class of biological systems is the subject
of an intensive area of research since a few years. Spatially extended
conductance based models of biological excitable membranes, and more
particularly of neuronal membranes \cite{HH52}, are considered, at an appropriate
scale, as hybrid stochastic models \cite{Aus08,FWL05}. More precisely, these models
can be regarded as belonging to the subclass of stochastic hybrid models
consisting of PDMPs \cite{BR11}.

The neurons are the basic building blocks of
the nervous systems. They are, roughly speaking, made of three parts: the
dendrites, the soma and the axon. In the present discussion, we focus our
attention on the axon, also referred as nerve fiber. The main role of the axon
is to transmit, for example to another neuron, the signal it received from the
soma. This electrical signal is called the nerve impulse or the action
potential or the spike. The target neuron may be far from the neuron which
emits the spike, therefore the nerve fiber is often much longer than width and
for mathematical convenience we model the axon by a segment, here the segment
$I=[0,1]$. As for all biological cells, the neuronal membrane allows the
exchanges of ions between the internal and the external cellular media through
the ionic channels. We assume that these ionic channels are disposed at
discrete location along the axon: at a location $z_i\in\xZ$ is a ionic channel
where $\xZ$ is a finite subset of $(0,1)$ of cardinality $N$. When an ionic
channel is open, it allows a flow of ions to enter or leave the cell, that is,
it allows a current to pass. We denote by $\xi\in E$, where $E$ is a finite
state space, the state of the channel at location $z_i$. Basically a state is
to be open or closed. When open, an ionic channel in the state $\xi$ allows the
ions of a certain specie associated to $\xi$ to pass through the membrane. If
we write $u(z_i)$ for the local potential of the neuronal membrane, often
called the trans-membrane potential, then the current which passes through the
membrane is given by
\[
  c_\xi(v_\xi-u(z_i))
\]
where $c_\xi$ is the conductance
of the ionic channel in state $\xi$ and $v_\xi$ the associated driven
potential. The quantities $c_\xi$ and $v_\xi$ are both real constants, the
first being positive. The sign of $v_\xi-u(z_i)$ tells us if the ions
corresponding to the state $\xi$ are leaving or entering the cell. The rate at
which the channel at location $z_i$ goes from one state $\xi$ to another state
$\zeta$ depends on the local trans-membrane potential $u(z_i)$ and is denoted
in the sequel by $q_{\xi\zeta}(u(z_i))$. For example, when the local potential
$u(z_i)$ is high, the ionic channels in a specific state $\xi$ will have the
tendency to open whereas when $u(z_i)$ is low, the ionic channels in the same
state $\xi$ will have the tendency to close. This is this mechanism which
allows the creation of an action potential and gradually its propagation along
the axon. Let $C$ be the capacitance of the axon membrane, $a$ the radius of
the axon and $R$ its internal resistance, this three quantities being positive
constants. The classical conductance based model for the propagation of the
nerve impulse is the cable equation
\begin{equation}\label{model:HH:pot}
C\partial_t u=\frac{a}{2R}\Delta u+\frac1N\sum_{i\in\xZ}c_{r(i)}(v_{r(i)}-u(z_i))\delta_{z_i}+\mathcal{I}
\end{equation}
where $r(i)$ denotes the state of the channel at location $z_i$. The dynamic of the channel $r(i)$ is given by
\begin{equation}\label{model:HH:ion}
  \prb{r_{t+h}(i)=\zeta|r_t(i)=\xi}=q_{\xi\zeta}(u_t(z_i))h+{\rm o}(h)
\end{equation}
where the ionic channels $\{r(i),i\in\xZ\}$ evolve independently over
infinitesimal time-scales. The quantity $\mathcal{I}$ in~\eqref{model:HH:pot}
corresponds to the impulse received by the axon. Since it does not affect our
analysis, we state $\mathcal{I}=0$ in the sequel. At last, the Partial
Differential Equation~\eqref{model:HH:pot} has to be endowed with boundary
conditions and we consider in the present paper zero Dirichlet boundary
conditions.

Let us be more specific about the assumptions on the jump rates
of the ionic channels. We assume that the function $q_{\xi\zeta}$ are smooth
and take values between two fixed positive constants. This assumption certainly
holds in the classical models for excitable membrane such as the Hodgkin-Huxley
model, Morris-Lecar model or models in Calcium dynamics. To gain a more
accurate insight in conductance based neuron model we refer to the classical
book~\cite{Hil84}.

Let us denote by $r=(r(i),i\in\xZ)$ the configuration of the whole ionic
channels: it is  a jump process with values in $E^N$. For a given
fixed trans-membrane potential~$u$, $r$ follows the following dynamic. For
$r\in E^N$, the total rate of leaving state~$r$ is given by
\[
\Lambda(u,r)=\sum_{i\in\xZ}\sum_{\xi\neq r(i)}q_{r(i),\xi}(u(z_i))
\]
and, if another state~$r'$ differs from $r$  only by the component $r(i_0)$,
then~$r'$ is reached with probability
\[
  Q(u,r)(\{r'\})=\frac{q_{r(i_0),r'(i_0)}(u(z_{i_0}))}{\Lambda(u,r)}
\]
If $r'$ differs from $r$ by two or more components then $Q(u,r)(\{r'\})=0$.
Between each jump of the process~$r$, the potential~$u$ follows the
deterministic dynamic induced by the PDE~\eqref{model:HH:pot}. For a fixed
time horizon~$T$, the process $((u_t,r_t),t\in[0,T])$ described this way
is a piecewise deterministic process (PDP) with values in the infinite dimensional
space\footnote{$H^1_0(I)$ denotes the usual Sobolev space of  functions in
$L^2(I)$ with first derivative in the sens of the distributions also in
$L^2(I)$ and with trace equals to zero on the boundary $\{0,1\}$ of $I$.}
$H^1_0(I)\times E^N$. This PDP is Markov --- and therefore a PDMP ---
if one chooses an appropriate filtration:
\begin{theorem}[\cite{BR11}]
There exists a filtered probability space satisfying the usual conditions such
that the process~$(u_t,r_t)$ is a homogeneous Markov process on~$H^1_0(I)\times E^N$,
\end{theorem}

\subsection{Singular perturbation of conductance based neuron model}

In this section we consider the spatially extended conductance based model of
the previous section as a slow-fast system: some states of the ion channels
communicate faster between them than others. This phenomenon is intrinsic to
classical neuron model, see \cite{Hil84}. Mathematically, this leads to
introduce a parameter $\e>0$, destined to be small, in equations
(\ref{model:HH:pot}) and (\ref{model:HH:ion}). For the states which communicate
at a faster rate, we say that they communicate at the usual rate divided by
$\e$.

The state space $E$ is therefore partitioned into different classes of states
which cluster the states communicating at a high rate. This kind of description
is usual, see for example \cite{FGC10}. We regroup our states in classes
making a partition of the state space $E$ into
\[
E=E_1\sqcup\cdots\sqcup E_l
\]
where $l\in\{1,2,\cdots\}$ is the number of classes. Inside a class $E_j$, the
states communicate faster at jump rates of order $\frac1\e$. States of
different classes communicate at the usual rate of order $1$. For $\e>0$ fixed,
we denote by $(u^\e,r^\e)$ the modification of the PDMP introduced in the
previous section with now two time scales. The new equations for the two
time-scales model are
\begin{equation}\label{model:HH:pot:e}
C\partial_t u^\e=\frac{a}{2R}\Delta u^\e+\frac1N\sum_{i\in\xZ}c_{r^\e(i)}(v_{r^\e(i)}-u^\e(z_i))\delta_{z_i}
\end{equation}
that is there is no changes in the form of the equation on the potential
whereas the dynamic of the channel $r^\e(i)$ is given by
\begin{equation}\label{model:HH:ion:e}
  \prb{r^\e_{t+h}(i)=\zeta|r^\e_t(i)=\xi}
  = q^\e_{\xi\zeta}(u^\e_t(z_i))h+{\rm o}(h)
\end{equation}
where $q^\e_{\xi\zeta}=\frac1\e q_{\xi\zeta}$ if $\xi$ and $\zeta$ are in the
same class and $q^\e_{\xi\zeta}=q_{\xi\zeta}$ otherwise.

For any fixed potential $u(z_i)=y$ we assume that the jump process $r(i)$
restricted to the class $E_j$ is irreducible and has therefore a unique
quasi-stationary distribution $\mu_j(y)$. The measure $\mu_j(y)$ does not
depend on $i$ because when the potential is held fixed to $y$, all the channels
$\{r(i),i\in\xZ\}$ have same law. The main idea is that when $\e$ goes to zero,
the fast components of the system (\ref{model:HH:pot:e})-(\ref{model:HH:ion:e})
reach their stationary behavior that can be viewed as an averaged dynamic. For
the jump process $r$ with  a potential held fixed to $y$, this leads to
consider the averaged jump process $\bar r=(\bar r(i),i\in\xZ)$ with values in
$\{1,\cdots,l\}^N$ and jump rates between two different classes given by
\[
\bar q_{l_1l_2}(y)=\sum_{\zeta\in E_{l_1}}\sum_{\xi\in E_{l_2}}q_{\zeta\xi}(y)\mu_{l_1}(y)(\zeta)
\]
$\bar q_{l_1l_2}(y)$ is indeed the sum of the rate of jumps from one state
$\zeta$ in $E_{l_1}$ to a state $\xi$ in $E_{l_2}$ averaged against the
quasi-stationary measure associated to the class $E_{l_1}$. The same idea can
be applied to the equation on the potential (\ref{model:HH:pot:e}) and
furthermore we have the following result.
\begin{theorem}[\cite{GT12}]
The process $(u^\e_t,t\in[0,T])$ converges in law in
$\mathcal{C}([0,T],H^1_0(I))$ when $\e$ goes to zero toward the process $u$
solution of
\begin{equation}\label{model:HH:pot:av}
C\partial_t u=\frac{a}{2R}\Delta u+\bar G_{\bar r}(u)
\end{equation}
with zero Dirichlet boundary conditions. The averaged reaction term is given by
\[
\frac{1}{N}\sum_{i\in \xZ}\sum_{j=1}^l 1_{j}(\bar r(i))\sum_{\xi\in E_j}\mu_j(u(z_i))(\xi)c_\xi(v_\xi-u(z_i))\delta_{z_i}
\]
and the dynamic of the averaged ionic channels is
\begin{equation}\label{model:HH:ion:av}
  \prb{\bar r_{t+h}(i)=l_2|\bar r_t(i)=l_1}=\bar q_{l_1 l_2}(u_t(z_i))h+{\rm o}(h)
\end{equation}
\end{theorem}
The above theorem can be seen as a law of large number as $\e$ tends
to zero. To continue our analysis, we study the fluctuations around
the averaged limit with a Central Limit Theorem (CLT).
\begin{theorem}
The fluctuations $z^\e=\frac{u-u^\e}{\sqrt\e}$ converge in law in
$\mathcal{C}([0,T],H^1_0(I))$ when $\e$ goes to zero toward a process $z$
corresponding to the solution of a Hybrid Stochastic Partial Differential
Equation (HSPDE).
\end{theorem}
One of the main advantages of considering the averaged model
(\ref{model:HH:pot:av})-(\ref{model:HH:ion:av}) instead of the two
time-scales model (\ref{model:HH:pot:e})-(\ref{model:HH:ion:e}) is that the
dynamic of the ionic channels is simplified: the state space for the averaged
ionic channels has a smaller cardinality than the original state space.
Averaging reduces in this case the dimension of the model. This can be very
useful for simulations or bifurcation analysis. However, we notice that the
averaged PDE on the potential (\ref{model:HH:pot:av}) may appear to be much
more non linear than in the non averaged case (\ref{model:HH:pot:e}). We also
remark that the error made by averaging is controlled by the CLT and that an
associated Langevin equation can be derived. The Langevin approximation
consists also in a HSPDE. If well posed and tractable, the Langevin
approximation may be regarded as a good compromise between the two-time scales
and the averaged model.

\section{Estimation of the distribution of the inter-jumping times}
\label{sec_estimation}

This section is dedicated to nonparametric statistics for piecewise
deterministic Markov processes. A suitable choice of the state space
and the main characteristics of the process provide stochastic models
covering a large number of applications, in reliability for instance
or in biology as mentioned before. In this context, it seems essential
to develop estimation methods for this class of processes.

In this part, we consider a piecewise deterministic Markov process~$\{X(t)\}$
defined on an open subset~$\mathcal{X}$ of a separable metric space
$(\mathcal{E},d)$. The conditional distribution of the inter-jumping
time $S_{i+1}$ given the previous post-jump location $Z_i$ is given by
\[
  \forall t\geq0,
  \quad
  \prb{S_{i+1}> t| Z_i}
  = \exp\left(  - \int_0^t \lambda(\phi(Z_i,s))\ud s\right)
  \ind{\{0\leq t<t^\star(Z_i)\}},
\]
where we recall that~$t^\star(x)$ is the deterministic exit time from
$\mathcal{X}$ starting from~$x$. One may associate with the compound function
$\lambda\circ\phi$ the conditional density~$f$ of the inter-jumping times,
given for any~$(x,t)\in \mathcal{X}\times\xR_+$ by:
\[
f(x,t) = \lambda(\phi(x,t))
\exp\left(  - \int_0^t \lambda(\phi(x,s))\ud s\right) .
\]
We investigate here the nonparametric
estimation of $f(x,t)$ for any $x\in \mathcal{X}$ and $0\leq t<t^\star(x)$ from
the observation of one trajectory of the process within a long
time. In the context of piecewise deterministic Markov processes, all
the randomness is contained in the embedded Markov chain
$(Z_n,S_n)$. As a consequence, our statistics may be computed from the
post-jump locations $Z_n$ and the inter-jumping times $S_n$. A
precise formulation of the results and the investigated methods may be
found in \cite{Az12b}. The general framework of metric state space has
been chosen in order to provide an estimate which could be computed in
applications involving state spaces of the form $\mathcal{X}=E\times M$,
where $M$ is continuous and $E$ is discrete.

This study is premised on a previous work of Azaïs and al. In
\cite{Az12a}, the authors investigate the nonparametric estimation of
the jump rate for a class of nonhomogeneous marked renewal
processes. This kind of stochastic models may be directly related with
piecewise deterministic Markov processes with constant flow. For this
class of processes, one may exhibit some continuous-time martingales,
which appear in the well-known framework of Aalen's multiplicative
intensity model (see \cite{AalPHD,Aal77,Aal78}). Nevertheless, this
idea is not relevant for studying the estimation of the jump rate for
piecewise deterministic Markov processes, since the conditional
independence $\esp{g(Z_{n+1})|Z_n} = \esp{g(Z_{n+1})|Z_n,S_{n+1}}$ is
not satisfied in the general case. Consequently, we propose to examine
the conditional distribution of $S_{n+1}$ given the previous post-jump
location $Z_n$ and the next post-jump location $Z_{n+1}$. We prove
that there exists, under some regularity assumptions, a mapping
$\widetilde{\lambda}$ from $\mathcal{X}\times \mathcal{X}\times\xR_+$
into $\xR_+$, such that for any $t\geq0$,
\[
\prb{S_{n+1}>t|Z_n,Z_{n+1}}
= \exp\left(
  - \int_0^t \widetilde{\lambda}(Z_n,Z_{n+1},s)\ud s\right)
\ind{\{0\leq t<t^\star(Z_n)\}} .
\]
In addition, we state that there exists a structure of continuous-time
martingale, which is similar to the one of the multiplicative
intensity model. Thanks to these considerations, we provide a way to
estimate an approximation $l$ of the jump rate $\widetilde{\lambda}$.

Our strategy for estimating the conditional probability density
function $f(x,t)$ consists in the introduction of a partition $(B_k)$
of the state space and two functions $l(A,B_k,t)$ and $H(A,B_k,t)$,
where $A$ is a set containing $x$. (The set $A$ is not allowed to
intersect the boundary of the state space $\mathcal{X}$.) They provide a way to
approximate the function of interest. Indeed, $f(x,t)$ is close to
\[
\sum_{k} l(A,B_k,t) H(A,B_k,t) ,
\]
if the partition $(B_k)$ and the set $A$ are thin enough. One may give
an interpretation of these two functions. First, $l(A,B_k,t)$ is an
approximation of the jump rate $\widetilde{\lambda}(x,y,t)$ of the
inter-jumping time $S_{n+1}$ given $Z_n=x$ and $Z_{n+1}=y$, under the
stationary regime, for $x\in A$ and $y\in B_k$. Roughly speaking,
$l(A,B_k,t)$ may be seen as the jump rate from $A$ to $B_k$ at time
$t$ under the stationary regime. Furthermore, $H(A,B_k,t)$ is exactly
the conditional probability $\prbPt{\nu}{S_{n+1}>t,Z_{n+1}\in B_k |
Z_n\in A}$, where $\nu$ denotes the invariant measure of the Markov
chain $(Z_n)$. We propose to estimate both the functions $l(A,B_k,t)$
and $H(A,B_k,t)$ in order to provide a consistent estimator of the
density of interest $f(x,t)$.

Firstly, our work is inspired by the smoothing methods proposed by
Ramlau-Hansen in \cite{Ram83} and the previous work of Azaïs et
al. \cite{Az12a}. We provide a nonparametric kernel estimator
$\widehat{l}_n(A,B_k,t)$ of $l(A,B_k,t)$. The keystone to state the
consistency of this estimator lies in the structure of continuous-time
martingale mentioned above. Furthermore, we estimate the conditional
probability $H(A,B_k,t)$ by its empirical version
\[
\widehat{H}_n(A,B_k,t)
= \frac{
          \sum_{i=0}^{n-1} \ind{\{Z_{i+1}\in B_k\}}\ind{\{S_{i+1}>t\}}
\ind{\{Z_i\in A\}}}
       {
          \sum_{i=0}^{n-1}\ind{\{Z_i\in A\}}}.
\]
The error between the function of interest and its estimate breaks
down into two parts: a deterministic error caused by the approximation
of $f(x,t)$ and a stochastic one in the estimation of this
approximation. Our major result of consistency is stated in the
following proposition.
\begin{prop}
For any $\varepsilon, \eta>0$, there exist an integer $N$, a set $A$
and a partition $(B_k)$ such that, for any $n\geq N$ and
$0<r_1<r_2<\inf_{x\in A} t^\star(x)$,
\[
\prb{\sup_{r_1\leq t\leq r_2}\abs{f(x,t)-\sum_k\widehat{l}_n(A,B_k,t)
\widehat{H}_n(A,B_k,t)}>\eta}<\varepsilon.
\]
\end{prop}

In addition, one may state the uniform convergence on every compact
subset $K\subset {\mathcal X}$. This nonparametric estimator of the conditional
distribution of the inter-jumping times is easy to compute from the
observation of the embedded chain $(Z_n,S_n)$ within a long
time. However, the choice of the set $A$ and the partition $(B_k)$
remains an open issue.

\section{Numerical methods.}
\label{sec_numerics}

A lot of numerical methods have been developed recently to simulate
diffusion processes and compute expectations, stopping times and other
interesting quantities (see references in \cite{dSDG10} for
example). But since PDMPs are in essence discontinuous at random
times, these results turn out to be too specific to be applied to
them. Besides, another important source of complication is the fact
that the transition semigroup $\{P_{t}\}_{t\in\mathbb{R}_{+}}$ of
$\{X(t)\}$ cannot be explicitly computed from the local
characteristics $(\phi, \lambda,Q)$ of the PDMP (see
\cite{CD99,CD08}). Therefore, it turns out to be hard to give an
explicit expression for the Markov kernel $P$ associated with the
Markov chain $\{X(t)\}$. Also, the Markov chain $\{X(t)\}$ is, in
general, not even a Feller chain (see \cite{Dav93}, pages 76 and 77).

On the other hand, PDMPs exhibit some nice specific properties. For
instance, all the randomness of the process can be described by the
discrete time (continuous state space) Markov chain $\{\Theta_{n}\}$
of the post jump locations and inter-jump times.

There exist essentially two families of numerical methods.  The
first one is based on the discretization the Markov kernel $Q$
\cite{CD88, CD89,Cos93}.  The second one relies on a discretization of
the Markov chain $\{\Theta_{n}\}$ using quantization \cite{pages98}.

Thanks to these methods one can solve optimal stopping problems
\cite{CD88,dSDG10,Gug86}, impulse control problems \cite{Cos93,CD89,
dSD12}, approximate distributions of exit times \cite{BdSD12} and
compute expectations of functionals of PDMP's \cite{BdSD12b}.


\bigskip
\noindent
\textbf{Acknowledgements.} The authors warmly thank Beno\^ite de
Saporta and Florent Malrieu for their kind and useful advices on the redaction
of these proceedings. They also thank Arnaud Guillin for his
enthusiastic and stimulating organization of the ``Journées MAS 2012''
in Clermont-Ferrand and of their subsequent proceedings.
The research of J.-B Bardet, A. G\'enadot, N. Krell and P.-A. Zitt is partly supported by the Agence Nationale de la
Recherche PIECE 12-JS01-0006-01. The research of Romain Azaïs was supported by ARPEGE program of the Agence Nationale de la Recherche, project FAUTOCOES, number ANR-09-SEGI-004.
\newcommand{\etalchar}[1]{$^{#1}$}
\providecommand{\bysame}{\leavevmode\hbox to3em{\hrulefill}\thinspace}
\providecommand{\MR}{\relax\ifhmode\unskip\space\fi MR }
\providecommand{\MRhref}[2]{%
  \href{http://www.ams.org/mathscinet-getitem?mr=#1}{#2}
}
\providecommand{\href}[2]{#2}

\end{document}